\documentclass[12pt]{article}
\usepackage{epsfig,amsmath,hyperref}
\usepackage[english]{babel}
\usepackage{color,colordvi}
\usepackage{amsfonts,amssymb,amsmath,amsopn,amsxtra}
\usepackage{mathrsfs}
\usepackage{mathtools}
\usepackage{enumerate}
\usepackage{bbm}
\usepackage{graphicx}
\usepackage{rotating}
\usepackage[normalem]{ulem}
\usepackage{slashed}

\DeclareMathOperator{\csch}{csch}

\newcommand{\R}{\mathbb{R}}

\newcommand{\e}{\mathrm{e}}
\newcommand{\OO}{\mathcal{O}}

\newtheorem{theorem}{Theorem}

\begin{document}

\vspace*{2em}

\begin{center}

\textbf{Ring chains with vertex coupling of a preferred orientation} \\[1.2em]
Marzieh Baradaran$^{1,*}$, Pavel Exner$^{2,3,**}$, Milo\v{s} Tater$^{2,3,***}$ \\[.5em]
\textit{$^1$Department of Mathematics, Faculty of Nuclear Sciences and Physical Engineering, \\ Czech Technical University, B\v rehov\'a 7, 11519 Prague, Czechia} \\[.3em]
\textit{$^2$Department of Theoretical Physics, Nuclear Physics Institute, Czech Academy \\ of Sciences, 25068 \v Re\v z near Prague, Czechia}\\[.3em]
\textit{$^3$Doppler Institute for Mathematical Physics and Applied Mathematics, \\ Czech Technical University, B\v rehov\'a 7, 11519 Prague, Czechia} \\[.3em]
\textit{$^*$marzie.baradaran@yahoo.com} \\
\textit{$^{**}$exner@ujf.cas.cz} \\
\textit{$^{***}$tater@ujf.cas.cz} \\[2em]
Received ...

\end{center}

\begin{quote}
We consider a family of Schr\"odinger operators supported by a periodic chain of loops connected either tightly or loosely through connecting links of the length $\ell>0$ with the vertex coupling which is non-invariant with respect to the time reversal. The spectral behavior of the model illustrates that the high-energy behavior of such vertices is determined by the vertex parity. The positive spectrum of the tightly connected chain covers the entire halfline while the one of the loose chain is dominated by gaps. In addition, there is a negative spectrum consisting of an infinitely degenerate eigenvalue in the former case, and of one or two absolutely continuous bands in the latter. Furthermore, we discuss the limit $\ell\to 0$ and show that while the spectrum converges as a set to that of the tight chain, as it should in view of a result by Berkolaiko, Latushkin, and Sukhtaiev, this limit is rather non-uniform.
\\[1em]
\end{quote}

\vspace{1.5em}

%%%%%%%%%%%%%%%%%%%%%%
\section{Introduction}
\label{s:intro}
\setcounter{equation}{0}

Quantum graphs came to recognition only slowly. First proposed by Linus Pauling \cite{Pa36} as a model for aromatic hydrocarbon molecules, they were briefly investigated \cite{RSch53} and then forgotten for more than three decades. Even when rediscovered in the late 1980s in connection with the progress in fabrication techniques that allowed to produce graphlike structures of semiconductors and other materials, they looked like a topic of limited importance. The development of the next three decades revealed many facets of the quantum graph theory including some deep mathematical questions going far beyond the original purpose; we refer to the recent monograph \cite{BK13} for a survey and a rich bibliography.

A characteristic property of quantum graphs is the non-uniqueness of their vertex coupling. The Schr\"odinger operators supported by graphs play the role of quantum mechanical Hamiltonians and as such they have to be self-adjoint. This means that the functions from their domain must satisfy at the vertices the conditions matching their values and first derivatives: in a vertex $v$ connecting $n$ graph edges they are generally of the form
 % ------------- %
\begin{equation}\label{genbc}
 (U-I)\psi(v)+i(U+I)\psi'(v)=0,
\end{equation}
 % ------------- %
where $U$ is an $n\times n$ unitary matrix. A natural question then arises about the meaning of such a coupling because each choice of the matrix $U$ defines a different physics. The most often considered case is the \emph{Kirchhoff coupling}\footnote{The name is unfortunate having in mind that Kirhhoff law in electricity means the current conservation and in the quantum case the conservation of probability current is equivalent to the self-adjointness. Several alternatives have been proposed but none of them stuck.} with the functions continuous and the sum of the derivatives vanishing, or the more general $\delta$ coupling \cite{Ex96} corresponding to $U= {2\over n+i\alpha}\mathcal{J}-I$, where $\mathcal{J}$ is the matrix whose all entries are equal to one and $\alpha\in\R$. One way to understand the meaning of a vertex coupling was proposed in early days \cite{RSch53}, namely to consider the dynamics in a family of thin tubes built around the graph `skeleton' and to look what happens it their widths shrink to zero. If the tubes boundary is Neumann, such a limit leads to the Kirchhof coupling \cite{Po11, RSch53}, the $\delta$ coupling is obtained by adding a properly scaled potential in the vertex region.

This is, however, only a small subset of the family described by the conditions \eqref{genbc}. The tube approximation can be worked out for \emph{any} of those couplings but the scheme has to be modified in a way which is far from being simple \cite{EP13}. First of all, the graph topology has to be changed locally by disconnecting the edges and connecting them pairwise by links whose lengths will go to zero in the limit. Furthermore, one has to add potentials, local scalar ones at the vicinity of the vertices and in the midpoints of the connecting links, and vector ones along each of the latter. They all are chosen to scale with respect to the tubes diameter in a particular nonlinear way adjusted to the coupling we need to approximate in the `shrinking limit'. This solves the question from the existence point of view, but the sketched, rather baroque construction is of a rather limited use from the practical point of view. An alternative is to adopt a pragmatic approach and to ask whether a particular vertex coupling could be useful in a specific physical model.

The motivation to analyze the coupling considered in this paper comes from a recent attempt to use quantum graphs to model the \emph{anomalous Hall effect} \cite{SK15}, that is, the occurrence of a voltage perpendicular to the current though the sample \emph{without the presence of a magnetic field}. While the mechanism of the usual Hall effect is well understood in both the classical and quantum cases, the anomalous effect is much less clear; it is conjectured that it comes from internal magnetization in combination with the spin-orbit interaction. The model proposed in \cite{SK15} regards the material supporting the 2D electron gas as `carpet' of rings connected mutually by a $\delta$ coupling at the points where they touch. Simple and elegant as it may be, however, it has a drawback: to obtain a nontrivial Hall conductance, the authors have to assume that the electron motion on the rings has a preferred direction which is something difficult to justify from the first principles. On the other hand, one can find among the conditions \eqref{genbc} couplings with such a property. As the simplest example, two of us studied in \cite{ET18} the matching referring to the matrix
 % ------------- %
$$
U= {\scriptsize \left( \begin{array}{ccccccc}
0 & 1 & 0 & 0 & \cdots & 0 & 0 \\ 0 & 0 & 1 & 0 & \cdots & 0 & 0 \\ 0 & 0 & 0 & 1 & \cdots & 0 & 0 \\ \cdots & \cdots & \cdots & \cdots & \cdots & \cdots & \cdots \\ 0 & 0 & 0 & 0 & \cdots & 0 & 1 \\ 1 & 0 & 0 & 0 & \cdots & 0 & 0
\end{array} \right)},
$$
 % ------------- %
which at a fixed momentum value, conventionally set to $k=1$, exhibits a `maximum' rotation; written in the components, the condition reads
 % ------------- %
\begin{equation}\label{coup}
(\psi_{j+1}-\psi_{j})+i(\psi_{j+1}^{\prime}+\psi_{j}^{\prime})=0,\quad j=1,\dots,n,
\end{equation}
 % ------------- %
and it is obvious that it is not invariant with respect to the time reversal operation represented in quantum mechanics by complex conjugation. The discussion of this coupling revealed interesting properties, in particular, it was found that the high-energy behavior of the on-shell scattering matrix, $S(k) = \frac{k-1 +(k+1)U}{k+1 +(k-1)U}$, depends on the parity $n$ of the vertex involved. This has spectral consequences which we have illustrated in \cite{ET18} using the examples of square and hexagonal lattice graphs.

The first aim of the present paper is to examine this effect of graphs periodic in a single direction. We are going to discuss graphs in the form of a periodic chain of loops connected either tightly or loosely through connecting links of the length $\ell>0$ assuming the coupling \eqref{coup} at all the graphs vertices. Since the graphs in question are periodic, the spectral analysis can be performed using the Floquet method \cite[Chap.~4]{BK13} writing the Hamiltonian in question as
 % -------------- %
\begin{equation}\label{Floquet}
H_\ell = \int_{-\pi}^\pi H_\ell(\theta)\, \mathrm{d}\theta
\end{equation}
 % -------------- %
where the fiber operator $H_\ell(\theta)$ acts on $L^2(C_\ell)$, where $C_\ell$ is the period cell and $C_\ell^*=[-\pi,\pi)$ is the dual cell, or Brillouin zone. Each of them has a purely discrete spectrum and the spectrum of $H_\ell$ then coincides with the union $\bigcup_{\theta \in C_\ell^*} \sigma(H_\ell(\theta))$. It is well known that the unique continuation principle is in general not valid in quantum graphs \cite[Sec.~3.4]{BK13} and it will be indeed the case, $\sigma(H_\ell)$ consists here of an absolutely continuous part and infinitely degenerate eigenvalues, or `flat bands' in the physicist terminology.

The indicated property of the vertex coupling \eqref{coup} is clearly manifested. While the positive spectrum of the tightly connected chain covers the entire halfline, the one of the loose chain is dominated by gaps as can be seen from the fact that the probability of being in the spectrum, in the sense of \cite{BB13} defined properly by relation \eqref{prob,l>0} below, \emph{is zero for any $\ell>0$}. This brings us to another interesting point. In a recent paper conditions have been found under which the spectrum of a quantum graph behaves `continuously' when the lengths of some of its edges tend to zero \cite{BLS19}. Our graphs satisfy the sufficient condition derived there and indeed, the spectrum of $H_\ell$ converges in the set sense to that of $H_0$, the Hamiltonian of the tightly connected chain, as $\ell\to 0$. Our result shows, however, that this convergence may be highly non-uniform.

%%%%%%%%%%%%%%%%%%%%%%%%%%%%%%%%%
\section{Tightly connected rings}
\label{s:tight}
\setcounter{equation}{0}

\begin{figure}[h]
\centering
\includegraphics[scale=.7]{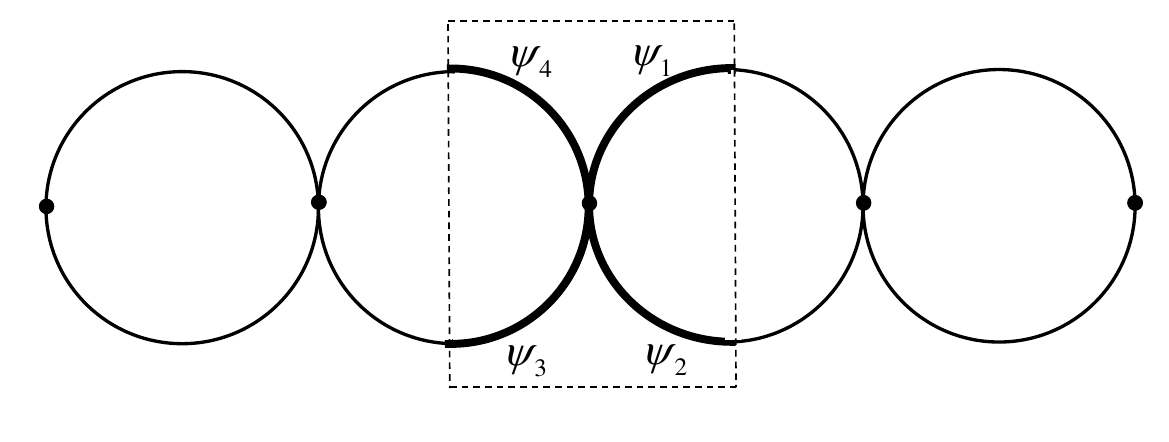}
\caption{ An elementary cell of the tightly connected ring chain}
\label{Fig1}
\end{figure}

Consider first the chain graph $\Gamma_0$ sketched in Fig.\ref{Fig1} in which the adjacent rings are coupled directly through the matching conditions \eqref{coup}; the corresponding Laplacian on $\Gamma_0$ will be denoted as $H_0$.
 % ------------- %
 % ------------- %
Following the standard Floquet-Bloch procedure \cite[Sec.~4.2]{BK13} we take an elementary cell of $\Gamma_0$ and  choose the coordinates to increase from left to right. Seeking a solution at energy $k^2$, we employ the following Ansatz
 % ------------- %
 \begin{align}\label{ansatz,l=0}
\psi_{j}(x)&=c_{j}^{+}e^{ikx}+c_{j}^{-}e^{-ikx},\quad\; x\in[0,\pi/2],\quad\;\;\: j=1,2, \nonumber \\[-.75em] \\[-.75em]
\psi_{j}(x)&=c_{j}^{+}e^{ikx}+c_{j}^{-}e^{-ikx},\quad\; x\in[-\pi/2,0],\quad j=3,4. \nonumber
\end{align}
 % ------------- %
The Floquet condition at the `free' ends of the cell gives the following conditions
 % ------------- %
\begin{equation}\label{Floq,l=0}
\psi_{j}(\pi/2)=e^{i\theta}\psi_{5-j}(-\pi/2), \quad  \psi^{\prime}_{j}(\pi/2)=e^{i\theta}\psi^{\prime}_{5-j}(-\pi/2), \quad j=1,2,
\end{equation}
 % ------------- %
with $\theta$ running through the dual cell, or Brillouin zone, $[-\pi,\pi)$. Substituting now (\ref{ansatz,l=0}) into \eqref{Floq,l=0}, one expresses coefficients $c_{j}^{\pm}$, $j=1,2$, in terms of $c_{j}^{\pm}$, $j=3,4$, as follows
 % ------------- %
\begin{equation} \label{onehalf}
c_{1}^{\pm}=c_{4}^{\pm}e^{i\theta}e^{\mp ik\pi},  \quad
c_{2}^{\pm}=c_{3}^{\pm}e^{i\theta}e^{\mp ik\pi}.
\end{equation}
 % ------------- %
Imposing next the condition (\ref{coup}) at the vertex, and taking into account that the derivatives in it are taken in the outward direction, we get
 % ------------- %
\begin{eqnarray}\label{sys,l=0}
\psi_{2}(0)-\psi_{1}(0)+i(\psi^{\prime}_{2}(0)+\psi^{\prime}_{1}(0))=0,\nonumber \\
\psi_{3}(0)-\psi_{2}(0)+i(-\psi^{\prime}_{3}(0)+\psi^{\prime}_{2}(0))=0, \nonumber \\[-.8em] \\[-.8em]
\psi_{4}(0)-\psi_{3}(0)+i(-\psi^{\prime}_{4}(0)-\psi^{\prime}_{3}(0))=0, \nonumber \\
\psi_{1}(0)-\psi_{4}(0)+i(+\psi^{\prime}_{1}(0)-\psi^{\prime}_{4}(0))=0. \nonumber
\end{eqnarray}
 % ------------- %
Inserting \eqref{ansatz,l=0} into\eqref{sys,l=0}, and taking into account \eqref{onehalf}, we get a system of four linear equations for the coefficients $c_{j}^{\pm}$, $j=3,4$. To be solvable, its determinant has to vanish, after an easy simplification it yields the spectral condition in the form
 % ------------- %
\begin{equation}\label{Sp.Con.P+,l=0}
k^3 \left(k^2+1\right) \sin k\pi\, \big(\cos k\pi - \cos\theta\big)=0.
\end{equation}
 % ------------- %
Furthermore, we know from \cite{ET18} that quantum graph operators with the coupling \eqref{coup} may not be positive, hence we have to inspect solutions with the energy $-\kappa^2$. This amounts to replacing a real $k$ by $k=i\kappa$ with $\kappa>0$, and the trigonometric functions in \eqref{Floq,l=0} by hyperbolic ones. In a similar way as above we then arrive at the spectral condition
 % ------------- %
\begin{equation}\label{Sp.Con.N-.l=0}
\kappa^3 \left(\kappa^2-1\right) \sinh\kappa\pi\, \big(\cosh\kappa\pi - \cos\theta\big)=0.
\end{equation}
 % ------------- %

 % ------------- %
\begin{theorem} \label{th:tight}
The spectrum of $H_0$ consists of the absolutely continuous part which coincides with the interval $[0,\infty)$, and a family of infinitely degenerate eigenvalues, the isolated one equal to $-1$ and the embedded ones equal to the positive integers.
\end{theorem}
 % ------------- %
\emph{Proof:} The condition \eqref{Sp.Con.N-.l=0} has the only solution, $\kappa=1$, while \eqref{Sp.Con.P+,l=0} is satisfied for all positive integers. The corresponding eigenvalues of the fiber operator are independent of $\theta$, hence the respective eigenvalues of $H_0$ have infinite multiplicity. The remaining part of the spectrum comes from the last factor at the left-hand-side of \eqref{Sp.Con.P+,l=0}. Every non-negative $k$ is obviously a solution corresponding to the energy $k^2 = \Big(\frac{\theta}{\pi}+2n\Big)^{2}$ with $n\in\mathbb{Z}$, and since the dependence on the quasimomentum $\theta$ is expressed by functions analytic on the Brillouin zone, this part of the spectrum is absolutely continuous. \hfill $\Box$

%%%%%%%%%%%%%%%%%%%%%%%%%%%%%%%%%
\section{Loosely connected rings}
\label{s:loose}
\setcounter{equation}{0}

Let us now replace $\Gamma_0$ by the graph $\Gamma_\ell$ in which the direct coupling of adjacent rings is replaced by connecting segments of length $\ell>0$ and the matching conditions \eqref{coup} are imposed at their ends; the corresponding Laplacian on $\Gamma_\ell$ will be denoted as $H_\ell$.
 % ------------- %
\begin{figure}[h]
\centering
\includegraphics[scale=.7]{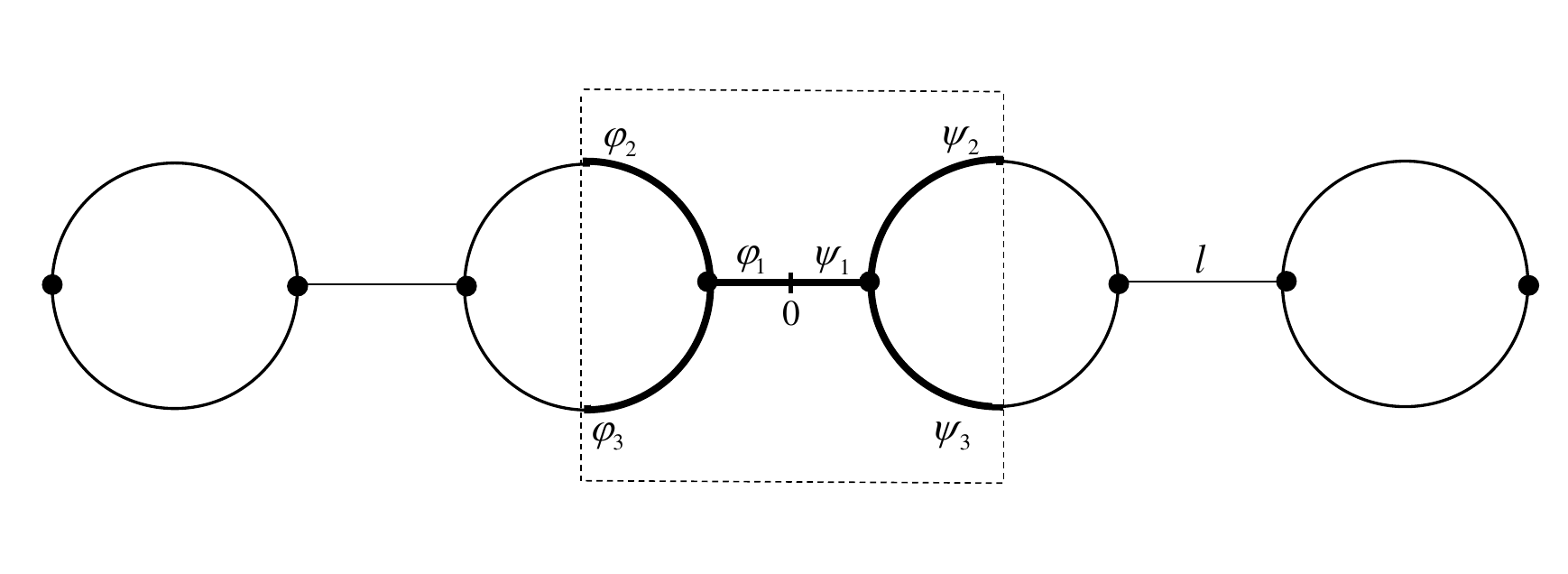}
\caption{The ring chain graph with $\ell>0$.}
\label{fig2}
\end{figure}
 % ------------- %
The Ansatz for positive energy solutions will in this case look as follows,
 % ------------- %
\begin{align}\label{ansatz,l>0}
&\psi_{1}(x)=a_{1}^{+}e^{ikx}+a_{1}^{-}e^{-ikx}, \quad\; x\in[0,\ell/2], \nonumber \\
&\psi_{j}(x)=a_{j}^{+}e^{ikx}+a_{j}^{-}e^{-ikx}, \quad\; x\in[0,\pi/2], \quad j=2,3, \nonumber \\[-.75em] \\[-.75em]
&\varphi_{1}(x)=b_{1}^{+}e^{ikx}+b_{1}^{-}e^{-ikx}, \quad\; x\in[-\ell/2,0], \nonumber \\
&\varphi_{j}(x)=b_{j}^{+}e^{ikx}+b_{j}^{-}e^{-ikx}, \quad\; x\in[-\pi/2,0], \quad j=2,3. \nonumber
\end{align}
 % ------------- %
Naturally, the functions $\psi_{1}$ and $\varphi_{1}$ have to be matched smoothly at segment midpoint, $\psi_{1}(0)=\varphi_{1}(0)$ and $\psi^{\prime}_{1}(0)=\varphi^{\prime}_{1}(0)$. On the other hand, the Floquet conditions at the boundary of the elementary cell require
 % ------------- %
\begin{equation}\label{Floq,l>0}
\psi_{j}(\pi/2)=e^{i\theta}\varphi_{j}(-\pi/2), \quad \psi^{\prime}_{j}(\pi/2)=e^{i\theta}\varphi^{\prime}_{j}(-\pi/2), \quad j=2,3.
\end{equation}
 % ------------- %
Using the smoothness at the midpoint together with \eqref{Floq,l>0}, we find
 % ------------- %
\begin{equation}\label{aij}
a_{1}^{\pm}=b_{1}^{\pm}, \quad
a_{j}^{\pm}=b_{j}^{\pm}e^{i\theta}\,e^{\mp ik\pi}, \quad j=2,3.
\end{equation}
 % ------------- %
The elementary cell now contains two vertices; the matching conditions (\ref{coup}) at them read
 % ------------- %
\begin{align}\label{sys,l>0}
&\psi_{3}(0)-\psi_{1}(\ell/2)+i(\psi^{\prime}_{3}(0)-\psi^{\prime}_{1}(\ell/2))=0, \nonumber \\
&\psi_{2}(0)-\psi_{3}(0)+i(\psi^{\prime}_{2}(0)+\psi^{\prime}_{3}(0))=0, \nonumber \\
&\psi_{1}(\ell/2)-\psi_{2}(0)+i(-\psi^{\prime}_{1}(\ell/2)+\psi^{\prime}_{2}(0))=0, \nonumber \\[-.75em] \\[-.75em]
&\varphi_{2}(0)-\varphi_{1}(-\ell/2)+i(-\varphi^{\prime}_{2}(0)+\varphi^{\prime}_{1}(-\ell/2))=0, \nonumber \\
&\varphi_{3}(0)-\varphi_{2}(0)+i(-\varphi^{\prime}_{3}(0)-\varphi^{\prime}_{2}(0))=0,\nonumber \\
&\varphi_{1}(-\ell/2)-\varphi_{3}(0)+i(\varphi^{\prime}_{1}(-\ell/2)-\varphi^{\prime}_{3}(0))=0. \nonumber
\end{align}
 % ------------- %
Plugging \eqref{ansatz,l>0} into \eqref{sys,l>0} and using \eqref{aij} we obtain a system of six linear equations; after a simple calculation one finds that it is solvable provided that the following condition is valid,
 % ------------- %
\begin{equation}\label{Sp.Con.l>0,Sin}
k^5 \sin k\pi \left(\left(k^4+2
   k^2+5\right) \sin k\pi \sin k\ell -4 \left(k^2+1\right)
   (\cos k\pi \cos k\ell -\cos \theta)\right)=0.
\end{equation}
 % ------------- %
As in the previous case the negative spectrum of $H_\ell$ corresponds to $k=i\kappa$ with $\kappa>0$; the solvability condition now acquires the form
 % ------------- %
 \begin{equation}\label{Neg,Sp,l>0}
\kappa ^5 \sinh \kappa\pi \left(4 \left(1-\kappa ^2\right) (\cosh \kappa\pi \cosh \kappa\ell- \cos \theta)+\left(\kappa ^4-2 \kappa ^2+5\right) \sinh \kappa\pi \sinh \kappa\ell \right)=0.
 \end{equation}
 % ------------- %
To state results, let us recall the notion of the \emph{probability of belonging the (positive) spectrum} put forward by Band and Berkolaiko \cite{BB13}: it is defined as
 % ------------- %
\begin{equation}\label{prob,l>0}
P_{\sigma}(H_\ell):=\lim_{K\to\infty} \frac{1}{K}\left|\sigma(H_\ell)\cap[0,K]\right|
\end{equation}
 % ------------- %
and it is clear that zero can be replaced here by any fixed positive number.
 % ------------- %
\begin{theorem} \label{th:loose}
The spectrum of $H_\ell$ has for any fixed $\ell>0$ the following properties:
 % ------------- %
\begin{enumerate}[(i)]
\setlength{\itemsep}{-3pt}
\item Any non-negative integer is an eigenvalue of infinite multiplicity.
\item Away of the non-negative integers the spectrum is absolutely continuous having a band-and-gap structure.
\item The negative spectrum is contained in $(-\infty,-1)$. It consists of a single band if $\ell=\pi$, otherwise there is a pair of bands and $-3\not\in {\sigma}(H_\ell)$.
\item The positive spectrum has infinitely many gaps.
\item $P_{\sigma}(H_\ell)=0$ holds for any $\ell>0$.
\end{enumerate}
 % ------------- %
\end{theorem}
 % ------------- %
\emph{Proof:} Everything follows, of course, from the conditions \eqref{Sp.Con.l>0,Sin} and \eqref{Neg,Sp,l>0}. The first two claims are easy. It is obvious from \eqref{Sp.Con.l>0,Sin} that non-negative integers solve the condition independently of $\theta$; in contrast to the tightly connected case there is no negative eigenvalue. The rest of the spectrum comes from the last factors at the left-hand sides of the said conditions; it is absolutely continuous since one can check, using the implicit function theorem, that the corresponding solutions depend analytically on the quasimomentum $\theta$.

Concerning the negative spectrum one may rewrite the condition \eqref{Neg,Sp,l>0} as
 % ------------- %
 \begin{equation}\label{Neg,Sp,l>0,cos}
\cos\theta=\cosh\kappa\ell\, \cosh\kappa\pi + \sinh\kappa\ell\, \sinh\kappa\pi\,\frac{\kappa^4-2 \kappa^2+5 }{4(1-\kappa^2)}.
\end{equation}
 % ------------- %
If $\kappa<1$, the last fraction is positive, and consequently, the right-hand side is larger than one, hence the interval $(-1,0)$ cannot belong to the spectrum. For $\kappa>1$ we can rewrite the condition \eqref{Neg,Sp,l>0,cos} as
 % ------------- %
 \begin{equation}\label{Neg,Sp,l>0,cos'}
\cos\theta= f_\ell(\kappa) := \cosh\kappa(\pi-\ell) - \frac{(\kappa^2-3)^2}{4(\kappa^2-1)} \sinh\kappa\ell\, \sinh\kappa\pi.
\end{equation}
 % ------------- %
It is obvious that
 % ------------- %
 \begin{equation}\label{l_extrem}
\lim_{\kappa\to 1+} f_\ell(\kappa) = \lim_{\kappa\to\infty} f_\ell(\kappa) = -\infty
\end{equation}
 % ------------- %
and that $f_\pi(\kappa) \le f_\pi(\sqrt{3}) =1$. To prove that for $\ell=\pi$ there is a single negative band we write $f_\pi(\kappa) = 1-h(\kappa)^2$ with $h(\kappa):= \frac12(\kappa^2-3)(\kappa^2-1)^{-1/2}\sinh \kappa\pi$ so that $h(\sqrt{3})=0$. The derivative of this function is
 % ------------- %
 \begin{equation}\label{h'}
h'(\kappa) = \frac{\kappa^2-3}{2\sqrt{\kappa^2-1}} \,\pi\, \cosh\kappa\pi + \frac{\kappa(\kappa^2+1)}{2(\kappa^2-1)^{3/2}}\, \sinh\kappa\pi,
 \end{equation}
 % ------------- %
hence we have $h'(\kappa)>0$ for $\kappa\ge\sqrt{3}$. On the other hand, $\lim_{\kappa\to 1+} h'(\kappa)=+\infty$, and consequently, $h(\cdot)$ would be also increasing in $(1,\sqrt{3})$ if there is no $\kappa$ in this interval for which one would have $h'(\kappa)=0$. Was it the case, the relation \eqref{h'} would lead to
 % ------------- %
 $$
\tanh\kappa\pi = \frac{\pi}{\kappa}\, \frac{(\kappa^2-1)(3-\kappa^2)}{\kappa^2+1}
 $$
 % ------------- %
which is impossible, since the right-hand side reaches its maximum value $\approx 0.812$ in the given interval at $\kappa\approx 1.303$ while $\tanh\kappa\pi > \tanh\pi \approx 0.996$.

On the other hand, for any positive $\ell\ne\pi$ we have $f_\ell(\sqrt{3}) = \cosh\kappa(\pi-\ell) > 1$, and consequently $-3\not\in \sigma(H_\ell)$, while in view of \eqref{l_extrem} and the continuity of $f_\ell(\cdot)$ the spectrum has a nonempty intersection with both the intervals $(-3,-1)$ and $(-\infty,-3)$. In each band the quasimomentum runs through the whole Brillouin zone and the band is determined by a pair of band edges; using \eqref{Neg,Sp,l>0,cos'} one can rewrite the conditions determining them, $f_\ell(\kappa)=\pm 1$, as
 % ------------- %
 \begin{equation}\label{alt}
\frac{(\kappa^2-3)^2}{4(\kappa^2-1)} = \frac{\cosh\kappa(\pi-\ell)\mp 1}{\sinh\kappa\ell\, \sinh\kappa\pi} =  \coth (\kappa \pi ) \coth (\kappa \ell )\mp\text{csch}(\kappa \pi ) \text{csch}(\kappa \ell )-1.
 \end{equation}
 % ------------- %
The left-hand side is positive, decreasing in $(1,\sqrt{3})$ and increasing in $(\sqrt{3},\infty)$, reaching its minimum at $\kappa=\sqrt{3}$ and diverging when $\kappa\to 1+$ or $\kappa\to\infty$. Since $\kappa \mapsto \coth\kappa\ell$ and $\kappa \mapsto \csch\kappa\ell$ are decreasing functions for any $\ell>0$, the right-hand side of \eqref{alt} with the lower sign is decreasing as a function of $\kappa$. As a result, the condition $f_\ell(\kappa)=-1$ has a single solution in $(\sqrt{3},\infty)$, and consequently, $H_\ell$ has a single band in the interval $(-\infty,-3)$.

In the interval $(1,\sqrt{3})$ this simple argument does not work; we shall instead show that the graph of the function $f_\ell$ can cross each of the values $\pm 1$ only once. To this aim we express the derivative of this function with respect to $\kappa$,
 % ------------- %
\begin{align*}
f_{\ell}^{\prime}(\kappa)&=-\frac{\left(\kappa ^2-3\right)^2}{4 \left(\kappa ^2-1\right)} \big( \ell  \sinh\kappa\pi\, \cosh\kappa\ell +\pi  \cosh\kappa\pi\, \sinh\kappa\ell\big)\\ \nonumber
&+\bigg(\frac{\kappa  \left(\kappa ^2-3\right)^2}{2 \left(\kappa ^2-1\right)^2}-\frac{\kappa  \left(\kappa ^2-3\right)}{\kappa ^2-1}\bigg) \sinh\kappa\pi\, \sinh\kappa\ell +(\pi -\ell ) \sinh \kappa(\pi-\ell ),
\end{align*}
 % ------------- %
which can be using simple manipulations brought to the form
 % ------------- %
\begin{align} \label{fp}
f_{\ell}^{\prime}(\kappa)&=-\frac{\left(\kappa ^2-3\right)^2}{4 \left(\kappa ^2-1\right)}\Big[\ell  \sinh\kappa\pi\, \cosh\kappa\ell+\pi  \cosh\kappa\pi\, \sinh\kappa\ell \nonumber \\ & +\frac{2 \kappa  \left(\kappa ^2+1\right)}{\left(\kappa ^2-3\right) \left(\kappa ^2-1\right)}\;\sinh\kappa\pi\, \sinh\kappa\ell\Big]+(\pi -\ell ) \sinh\kappa (\pi-\ell ).
\end{align}
 % ------------- %
Substituting now from (\ref{alt}) into (\ref{fp}) we find that
 % ------------- %
\begin{align}\label{fpsubp}
f_{\ell}^{\prime}(\kappa)&=\bigg(\frac{2 \kappa  \left(\kappa ^2+1\right)}{\left(3-\kappa ^2\right) \left(\kappa ^2-1\right)}-\pi  \coth\kappa\pi -\ell  \coth\kappa\ell\bigg)\big(\cosh \kappa(\pi-\ell )\mp1\big)\nonumber \\ & \quad +(\pi -\ell ) \sinh\kappa(\pi-\ell)
\end{align}
 % ------------- %
must hold if $f_{\ell}(\kappa)=\pm 1$, respectively. Our aim is to show that these derivatives are positive for any $\ell>0$ and $\kappa\in(1,\sqrt{3})$. Consider the case with the upper sign and suppose that
 % ------------- %
\begin{equation}\label{contr}
\Big[\frac{2 \kappa  \left(\kappa ^2+1\right)}{\left(3-\kappa ^2\right) \left(\kappa ^2-1\right)}-\pi  \coth\kappa\pi\, -\ell  \coth\kappa\ell )\Big]\big(\cosh\kappa (\pi-\ell )-1\big) \leq (\ell-\pi) \sinh\kappa(\pi-\ell )
\end{equation}
 % ------------- %
holds. Since $\cosh\kappa(\pi-\ell)-1>0$, we can rewrite this relation as
 % ------------- %
\begin{equation} \label{contr2}
\frac{2 \kappa  \left(\kappa ^2+1\right)}{\left(3-\kappa ^2\right) \left(\kappa ^2-1\right)}\leq\pi  \coth\kappa\pi +\ell  \coth\kappa\ell +(\ell-\pi) \coth \frac{\kappa(\pi-\ell)}{2}.
\end{equation}
 % ------------- %
We denote the left- and right-hand side of \eqref{contr2} as $g_1(\kappa)$ and $g_2(\kappa,\ell)$, respectively. They are both positive, for $g_2(\kappa)$ it is obvious for $\ell\ge\pi$, for smaller values the function is increasing as we will see in a minute, and $g_2(\kappa,0+)= \kappa^{-1} +\pi\big(\coth\kappa\pi - \coth\frac{\kappa\pi}{2}\big)$ is in the interval $(1,\sqrt{3})$ larger than $g_2(\sqrt{3},0+)\approx 0.550$. The function $g_1$ satisfies
 % ------------- %
$$
\lim_{\kappa\to 1+} g_{1}(\kappa)=\lim_{\kappa \to \sqrt{3} -} g_{1}(\kappa)=\infty
$$
 % ------------- %
and reaches its minimum value $\approx7.737$ at $\kappa\approx 1.303$.

On the other hand, asking about the values of $g_2(\kappa,\ell)$ we first note that for a fixed $\kappa\in(1,\sqrt{3})$ it is increasing as as function of $\ell$. To check that, consider the partial derivative
 % ------------- %
\begin{equation} \label{partd}
\frac{\partial g_{2}(\kappa,\ell)}{\partial \ell}=\frac{\kappa(\ell-\pi )}{\cosh\kappa(\pi-\ell)-1}+\coth\frac{\kappa(\pi-\ell)}{2}+\coth\kappa\ell- \kappa\ell\, \text{csch}^2\kappa \ell.
\end{equation}
 % ------------- %
The last two terms at the right-hand side in combination vanish in the limit $\ell\to 0+$ so for small $\ell$ the expression is dominated by the first two being thus positive. The same is true for all $\ell<\pi$, because was it not the case we would have
 % ------------- %
$$
\frac{\sinh \kappa(\pi-\ell)-\kappa(\pi-\ell)}{\cosh\kappa(\pi-\ell)-1} \leq  \kappa\ell\,\text{csch}^2 \kappa \ell -\coth\kappa\ell,
$$
 % ------------- %
or equivalently
 % ------------- %
$$
\frac{\sinh^2\kappa\ell}{\cosh\kappa(\pi-\ell)-1}\,\big[\sinh\kappa(\pi-\ell)-\kappa(\pi-\ell)\big] \leq -\big[\sinh\kappa \ell\, \cosh\kappa\ell -\kappa\ell\big],
$$
 % ------------- %
which is impossible because $\cosh\kappa(\pi-\ell)>1$ for $\ell<\pi$ and the expressions in the square brackets are obviously positive. Using further l'Hospital rule, we find that the derivative is also positive at $\ell=\pi$,
 % ------------- %
$$
\frac{\partial g_{2}(\kappa,\pi)}{\partial \ell} = \coth\kappa\pi- \kappa\pi\, \text{csch}^2\kappa \pi > \coth\pi- \pi\, \text{csch}^2\pi \approx 0.980
$$
 % ------------- %
To deal with the case $\ell>\pi$, we rewrite the condition \eqref{partd} as
 % ------------- %
\begin{equation} \label{partd2}
\frac{\partial g_{2}(\kappa,\ell)}{\partial \ell} = F(2\kappa\ell) - F(\kappa(\ell-\pi)),
\end{equation}
 % ------------- %
where
 % ------------- %
$$
F(u) := \frac{\sinh u - u}{2\sinh^2\frac{u}{2}}.
$$
 % ------------- %
It is easy to see that $F(0+)=0$ and $\lim_{u\to\infty} F(u)=1$ and that $F$ is increasing because
 % ------------- %
$$
F'(u) = \frac{u\,\sinh\frac{u}{2}-2}{\cosh u -1} > 0
$$
 % ------------- %
holds in view of the inequality $\coth x> x^{-1}$ valid for $x>0$. This monotonicity in combination with \eqref{partd2} shows that the derivative is positive for any $\ell>0$. It follows that
 % ------------- %
\begin{equation} \label{g_2}
g_2(\kappa,\ell) < \lim_{\ell\to \infty} g_{2}(\kappa,\ell) = \pi(1+\coth\kappa\pi) < \pi(1+\coth\kappa\pi) \approx 6.295
\end{equation}
 % ------------- %
Putting these results together we conclude that assumption \eqref{contr} leads to a contradiction, hence $f'_\ell(\kappa)>0$ holds at any point $\kappa\in(1,\sqrt{3})$ where $f_\ell(\kappa)=1$. In the same way we can check that the derivative is positive when the graph of the function crosses the value $-1$. This completes the proof of the part (iii); let us note that the same method can be used as an alternative way to prove the existence of a single band below $-3$ by verifying that in this case the derivatives $f'_\ell(\kappa)$ at the point referring to the band edges are negative.

Turning to claims \emph{(iv)} and \emph{(v)} about the positive spectrum we note first that for large values of $k$ the first term in the bracket at the left-hand side of the condition \eqref{Sp.Con.l>0,Sin} is dominating over the second, $\theta$-dependent one, and consequently, one can expect the spectral bands to become almost flat at high energies. To see that this is indeed the case and to give a more precise meaning to such an asymptotics,
we rewrite \eqref{Sp.Con.l>0,Sin} in the form
 % ------------- %
 \begin{equation}\label{Sp.Con.l>0,cot,alter1}
\cos k\ell \cos k\pi -\cos\theta = \frac{k^4+2 k^2+5}{4 \left(k^2+1\right)}\sin k\ell \sin k\pi.
\end{equation}
 % ------------- %
It shows that the points $k$ where $\sin k\ell \sin k\pi=0$ belong to the spectrum, since one can always find $\theta$ such that $\cos k\ell \cos k\pi =\cos\theta$. Asking about the left and right endpoints of the bands, we arrive at
 % ------------- %
  \begin{equation}
\cos k\ell \cos k\pi -1 \leq \frac{k^4+2 k^2+5}{4 \left(k^2+1\right)}\sin k\ell \sin k\pi \leq \cos k\ell \cos k \pi+1
\end{equation}
 % ------------- %
which implies that the middle expression takes values in the interval $[-2,2]$. By negation, a sufficient condition for being in the spectral gap is
 % ------------- %
 \begin{equation}\label{Gap.Con.l>0}
\frac{k^4+2 k^2+5}{4 \left(k^2+1\right)}\,|\sin k\ell \sin k\pi)|>2,
\end{equation}
 % ------------- %
which for can be replaced by a weaker sufficient condition with the same asymptotics for large $k$,
 % ------------- %
  \begin{equation}\label{Gap.Con.largek}
|\sin k\ell \sin k\pi |>\frac{8}{k^{2}}.
\end{equation}
 % ------------- %
Given now positive an interval $J=(k_{1},k_{2})\subset \mathbb{R}_+$ one can introduce the set
 % ------------- %
\begin{equation}\label{M_l}
M_{\ell}(k_{1},k_{2}):=\bigg\{k\in(k_{1},k_{2})\,:\,|\sin k\ell|>\frac{2\sqrt{2}}{k}\bigg\};
\end{equation}
 % ------------- %
we have clearly $k\not\in\sigma(H)$ if $k\in M_{\ell}(k_{1},k_{2})\cap M_{\pi}(k_{1},k_{2})$. Let us assess the measure of $M_{\ell}(k_{1},k_{2})$. We have $\sin(k\ell )=\pm\sin(\delta\ell )$ for $k=\frac{m\pi}{\ell}+\delta \in J$, and since $|\sin x|\geq\frac{2}{\pi}|x|$ holds for $|x|\le \frac12\pi$, the inequality $\delta>\frac{\pi\sqrt{2}}{k\ell}$ ensures that $k\in M_{\ell}(k_{1},k_{2})$. The interval $J$ contains a finite number of zeros, hence there is a constant $c>0$ depending on $|J|=|k_2-k_1|$ such that
 % ------------- %
\begin{equation}
|M_{\ell}(k_{1},k_{2})|\geq |k_{2}-k_{1}|-\frac{c}{k_1},
\end{equation}
 % ------------- %
and consequently
 % ------------- %
\begin{equation}|M_{\ell}(k_{1},k_{2})\cap M_{\pi}(k_{1},k_{2})|\geq |k_{2}-k_{1}|-\frac{2c}{k_1}.
\end{equation}
 % ------------- %
This allows us to find the probability $P_{\sigma}(H_\ell)$. We choose $k_{m}=md$ for $d>0$, then we have
 % ------------- %
\begin{align}
&1- \frac{1}{(N-1)d} \big| \sigma(H_\ell) \cap [d,Nd]\big| \geq \frac{1}{k_{N}-k_{1}}\sum_{m=1}^{N-1}\bigg(k_{m+1}-k_{m}-\frac{2c}{k_{m}}\bigg)\\\nonumber
& \qquad =\frac{1}{(N-1)d}\,\sum_{m=1}^{N-1}\bigg(d-\frac{2c}{md}\bigg)
=1-\frac{2c}{d^{2}}\,\frac{1}{N-1}\sum_{m=1}^{N-1}\frac{1}{m} \longrightarrow 1
\end{align}
 % ------------- %
as $N\to\infty$ because the sum equals $\gamma_\mathrm{E} + \psi(N) = \gamma_\mathrm{E} + \ln N + \mathcal{O}(N^{-1})$ where $\gamma_\mathrm{E}\approx 0.577$ is the Euler-Mascheroni constant; note that the constant $c$ is independent of $m$ because all the involved intervals have the same length. This yields the claim \emph{(v)} and at same time \emph{(iv)} because the operator $H_\ell$ is unbounded having all the positive integers in its spectrum; a finite number of gaps would then imply $P_{\sigma}(H_\ell)=1$. In this way, the proof is concluded. \hfill $\Box$

\medskip

The last claim of the theorem means that, in contrast to the case of tightly connected rings, the transport over such a chain is impossible for most energies because the spectrum is dominated by the gaps. This feature illustrates one more time the role of the vertex parity observed in \cite{ET18}. The asymptotic behavior of individual bands at high energies may differ in dependence on the parameters involved. If the right-hand side of \eqref{Sp.Con.l>0,cot,alter1} has a double root, the band width in the $k$ variable is proportional to $k^{-1}$, and consequently, in the energy variable it is of order $\mathcal{O}(1)$ as $k\to\infty$. This is analogous to the behavior of the $\delta'$ Kronig-Penney model \cite[Chap.~III.3]{AGHH} and its lattice generalizations \cite{Ex96}. Note that such bands have the infinitely degenerate eigenvalues from claim \emph{(i)} of Theorem~\ref{th:loose} embedded in them. On the other hand, if the root in \eqref{Sp.Con.l>0,cot,alter1} is simple -- which happens always, for instance, if $\frac{\ell}{\pi}\not\in\mathbb{Q}$ -- the corresponding band widths are of order $\mathcal{O}(k^{-1})$ as $k\to\infty$.

%%%%%%%%%%%%%%%%%%%%%%%%%%%%%%%%%%%%%%%%%%%%%%
\section{Changing the connecting links length}
\label{s:segment}
\setcounter{equation}{0}

Let us first consider the situation when the segments connecting the rings shrink to points, i.e. $\ell\to 0$. The spectral condition \eqref{Sp.Con.l>0,cot,alter1} turns in the limit into $\cos k\pi-\cos\theta=0$ which leads to the same conclusion as \eqref{Sp.Con.P+,l=0} suggesting that the limiting positive spectrum fills the interval $[0,\infty)$, naturally in addition to the infinitely degenerate eigenvalues at the positive integers which are $\ell$-independent. This conclusion is not surprising taking into account the results of \cite{BLS19}. Applying Theorem~3.5 of this work to the fiber operator in \eqref{Floquet} we see that $H_\ell(\theta)\to H_0(\theta)$ in the norm-resolvent sense, since the domain of $H_\ell$ contains no function that would be nonzero only on the connecting link, hence Condition~3.2 of \cite{BLS19} is satisfied. It follows that the eigenvalues of $H_\ell(\theta)$ converge to those of $H_0(\theta)$, and in view of \eqref{Floquet} the spectrum of $H_\ell$ converges to that of $H_0$ in the sense of sets.

At the same time, the convergence is \emph{highly nonuniform}. Indeed, by claim~\emph{(v)} of Theorem~\ref{th:loose} the probability \eqref{prob,l>0} that a positive energy belongs to the spectrum is \emph{zero} for any $\ell$ while $P_\sigma(H_0)=1$ holds by Theorem~\ref{th:tight}. The reason is clear: it follows from \eqref{Sp.Con.l>0,cot,alter1} that $k^2$ belongs to a gap if and only if
 % ------------- %
\begin{equation}\label{Gap.Con.l>0'}
\big|\cos k\ell \cos k\pi - \frac{k^4+2 k^2+5}{4 \left(k^2+1\right)}\sin k\ell \sin k\pi\big|>1.
\end{equation}
 % ------------- %
It is the second term on the left-hand side which is responsible for opening the gaps, and the smaller $\ell$ is the larger $k$ must be to make it dominate over the first one.

The behavior of the negative spectrum is equally interesting; what is left of it in the limit is the value $-1$ which is, as we know, the eigenvalue of the four-prong star graph with the coupling \eqref{coup}. Let us summarize the results:
 % ------------- %
\begin{theorem} \label{th:l->0}
In the limit $\ell\to 0$ the positive spectrum of $H_\ell$ fills the interval $[0,\infty)$. The lower edge of the upper negative band behaves asymptotically as
 % ------------- %
\begin{equation}\label{up-low-asy}
-\kappa^2(\ell;0) = -1 -\ell\, \coth\frac{\pi}{2} + \OO(\ell^2) \quad \text{for} \;\; \ell\to 0+
\end{equation}
 % ------------- %
and the band width is $-\kappa^2(\ell;\pm\pi)+\kappa^2(\ell;0) = 2\ell(\sinh\pi)^{-1} + \OO(\ell^2)$. On the other hand, the lower negative band escapes to $-\infty$, more specifically, we have
 % ------------- %
\begin{equation}\label{low-asy}
-\kappa^2(\ell;\theta) = -\Big(\frac{4}{\ell}\Big)^{2/3}\!\! + \OO(1) \quad \text{for} \;\; \ell\to 0+
\end{equation}
 % ------------- %
and any fixed quasimomentum $\theta\in[-\pi,\pi)$ with the error term independent of $\theta$.
\end{theorem}
 % ------------- %
\emph{Proof:} We have already dealt with the positive spectrum. Let us next rewrite the spectral condition \eqref{Neg,Sp,l>0,cos'} as an implicit equation 
 % ------------- %
 \begin{equation}\label{Neg,Sp,implicit}
0 = g(\kappa,\ell):= (\kappa^2-3)^2 \sinh\kappa\ell\, \sinh\kappa\pi + 4(\kappa^2-1)(\cos\theta - \cosh\kappa(\pi-\ell))
\end{equation}
 % ------------- %
We have $g(1,0)=0$ and
 % ------------- %
$$
\frac{\partial g}{\partial\kappa}(1,0) = 8(\cos\theta - \cosh\pi) \ne 0,
$$
 % ------------- %
which justifies the use of implicit function theorem. Since $\frac{\partial g}{\partial\ell}(1,0) = 4\sinh\pi$, it yields
 % ------------- %
 \begin{equation}\label{Neg,Sp,upper}
\kappa(\ell;\theta) = 1 + \frac{\ell}{2}\, \frac{\sinh\pi}{\cosh\pi-\cos\theta} + \OO(\ell^2) \quad \text{as} \;\; \ell\to 0+
\end{equation}
 % ------------- %
which implies the claims made about the upper negative band. To see the behavior of the lower one, let us first note that the corresponding function must satisfy $\lim_{\ell\to 0}\kappa(\ell;\theta)=\infty$, otherwise the condition \eqref{Neg,Sp,l>0,cos'} could not hold for $\ell$ small enough. Furthermore, we rewrite the condition \eqref{Neg,Sp,l>0,cos} as
 % ------------- %
 \begin{equation}\label{Neg,Sp,l>0,cos''}
\frac{\kappa^4-2 \kappa^2+5 }{4(1-\kappa^2)}\,\tanh\kappa\ell\, \tanh\kappa\pi = \frac{\cos \theta}{\cosh \kappa\ell \cosh \kappa\pi}-1.
\end{equation}
 % ------------- %
This implies $\lim_{\ell\to 0}\kappa(\ell;\theta)\ell=0$ because otherwise there would be a sequence $\{\ell_n\},\: \ell_n\to 0$, for which the condition will have no solution but we know that it is not the case. As a result, we have for small $\ell$ and large $\kappa$
 % ------------- %
 $$
\frac14 \kappa^2\cdot \kappa\ell\, \big(1+\OO(\kappa^{-2}) + \OO(\kappa\ell)\big) = 1 + \OO(\kappa^2\ell^2)
 $$
 % ------------- %
which yields
 % ------------- %
 \begin{equation}\label{Neg,Sp,upper}
\kappa(\ell;\theta) = \sqrt[3]{\frac{4}{\ell}} + \OO(\ell^{1/3}) \quad \text{as} \;\; \ell\to 0+
\end{equation}
 % ------------- %
and consequently, the relation \eqref{low-asy} independently of $\theta$. \hfill $\Box$

\medskip

Let us now turn to the opposite extremal situation. The limit $\ell\to\infty$ makes little sense but we can ask about the spectral behavior of the system when the spacing $\ell$ is large. While the gaps dominate according to Theorem~\ref{th:loose}, the spectrum becomes nevertheless `dense' in the sense that the gap with a fixed ordinal number can be made arbitrarily narrow by choosing $\ell$ large enough. Indeed, we know from the proof of the said theorem that the points with $k=\frac{\pi n}{\ell},\: n\in\mathbb{N}$, belong to the spectrum, and this implies that the distance of the neighboring bands is $\:\approx (2n+1)\big(\frac{\pi}{\ell}\big)^2 + \mathcal{O}(l^{-1})$ for large $\ell$.

The negative spectrum bands shrink in the asymptotic regime, however, in contrast to the limit $\ell\to 0+$ not to the same point. In order to see that, we rewrite the condition \eqref{Neg,Sp,l>0,cos'} in still another form,
 % ------------- %
 \begin{equation}\label{Neg,Sp,l>0,cos''}
\left(1+\frac{(\kappa^2-3)^2}{4(\kappa^2-1)}\right) \tanh\kappa\pi = \coth\kappa\ell - \frac{\cos\theta}{\sinh\kappa\ell\,\cosh\kappa\pi}.
\end{equation}
 % ------------- %
For large $\ell$ the right-hand side behaves as $1+\mathcal{O}(\e^{-\ell})$. To find the approximate band positions we put $\kappa^2=3+\varepsilon$, the the expression in the bracket on the left-hand side is $1+\frac{\varepsilon^2}{4(2+\varepsilon)}$, and putting $\tanh\kappa\pi \approx 1-\e^{-2\pi\sqrt{3}}$ we find that
 % ------------- %
 \begin{equation}\label{epsilon}
\varepsilon \approx \pm 4\,\e^{-\pi\sqrt{3}} \approx \pm  0.0173.
 \end{equation}
 % ------------- %
This is not surprising, recall that that the three-edge star graph with the vertex coupling \eqref{coup} has the eigenvalue $-3$ with the eigenfunction localized exponentially in the vicinity of the vertex. A similar squeezing of a band pair around this value occurs in hexagonal lattice graphs with long edges \cite{ET18}. Here, however, the two bands do not converge to a single point but a pair of separate ones, the reason being that the lengths of two of the three edges meeting in a vertex are fixed, albeit large in comparison with the `size' of the indicated eigenfunction, and only one is getting longer.

\subsection*{Acknowledgements}
The research was supported by the Czech Science Foundation (GA\v{C}R) within the project 17-01706S and by the EU project CZ.02.1.01/0.0/0.0/16\textunderscore 019/0000778. We are obliged to the referee for useful comments.

\end{document}